\documentclass{article}
\usepackage{authblk}
\usepackage[english]{babel}
\usepackage{amsthm}
\usepackage{amsmath}
\usepackage{amssymb}
\usepackage{amsfonts}
\usepackage[square,numbers]{natbib}
\usepackage{hyperref}
\usepackage{multirow}
\usepackage{babel}
\usepackage[T1]{fontenc}
\newtheorem{theorem}{Theorem}[section]

\newtheorem{lemma}[theorem]{Lemma}
\def\ff{\ensuremath{\mathbf{F}}}

\title{Construction of Permutation Polynomials of Certain 
Specific Cycle Structure over Finite Fields}

\author[1]{Anitha G}
\author[2]{P Vanchinathan}
\affil[1,2]{School of Advanced Sciences\\ 
 Vellore Institute of Technology\\ 
 Chennai 600 127\\ INDIA\\}
\affil[1]{\texttt{anitha.g2019@vitstudent.ac.in} }
\affil[2]{\texttt{ vanchinathan.p@vit.ac.in} }

\begin{document}
\date{}
\maketitle
\begin{abstract}
For a finite field of order an odd number $q$, 
	we construct families of permutation binomials and permutation trinomials
	with one fixed-point (namely zero) and
	remaining elements being permuted as  disjoint  cycles  of same length 
	(i.e. cycle type $1+m^d$).
	Further, a family of  permutations having cycle type of the form $1+m^r+n^s$ with $m,n,r,s$ all divisors of $q-1$
	 are constructed as binomial and trinomial permutations.
	 A generalization to polynomials
	 of more number of terms is also provided.
\end{abstract}

\section{Introduction}
Polynomials defined over finite fields which are bijective functions on that field
have been an object of study  by Hermite\cite{hermite} and Dickson \cite{dickson} in late nineteenth century. 
Called permutation polynomials they have found applications in cryptography, algebraic coding theory and in combinatorics.

Any automorphism of the multiplicative group extended to the whole field by sending 0 to 0, is
given by monomial $f(x) = x^k$, with $k$ coprime to the order of the multiplicative group.

These monomials along with the linear polynomials $ax+b,\ a\ne0$ are readily the easiest permutation polynomials one can think of.

More and more families of permutation polynomials have been found by various mathematicians. Chapter 8 of the
monograph by Mullen and Panerio \cite{handbook}  and Chapter 7 of Lidl and Niedereiter \cite{lidl}  and the survey article \cite{hou} by Hou  provide a wealth of information
on varieties of permutation polynomials. Due to ease of applying to electronic hardware,
permutation polynomials for fields of characteristic 2 have received lot of attention. 
(See, for example, 
\cite{E1,E3,E5}.)

Here we focus on polynomials with fewer terms, and work in odd characteristic.
In real life applications, as the time needed to evaluate a polynomial is a crucial factor, it is desirable to find polynomials with fewer terms. Naturally there  have been attempts to
identify all permutation polynomials with two terms, the binomials.

Permutation binomials with terms involving just $X^{(q+1)/2}$  and $X$
over a finite field $\ff_q$ have been completely determined. (See page 352, Theorem 7.11 in \cite{lidl}).
A result of Carlitz and Wells \cite{carwells} assures existence of permutation binomials
of the form $x^{d+1} + ax$ for all  $d|q-1$.
More results on binomial permutations are in in Lappano \cite{lap},  Hou and 
Lappano \cite{houlap}, and \cite{nieder}.

Construction of permutation trinomials can be found in the works of  
Bartoli and Timponella \cite{bart},  Gupta and Sharma\cite{rohit} (in characteristic 2), 
Hou, \cite{hou1,hou2},  Lavorante\cite{lavo}, and Li, Qu and Chen \cite{lqc}.

However, the determination of the  cycle type as a permutation  for known permutation polynomials
is a difficult task. There does not seem to be any easily discernible
relationship between  polynomials on the one hand,  and
the cycle types of the permutations they provide on the other.
Even the problem of  giving a  family of polynomials all giving rise to permutations of
some specific cycle type is stated to be difficult [See \cite{redei} ] 
as recently as in the year 2022. 
Literature on this question is very limited. 
\textit{In this paper we report our findings on constructing families of permutation polynomials which are simultaneously simple as  polynomials (binomial, trinomial etc)   and have uncomplicated 
cycle type -- as partitions they use just 2 or 3 numbers repeatedly}.

\textit{At the outset we would like to mention that all the permutation 
polynomials constructed in this paper
fall under the category of cyclotomic mappings.}  (See \cite{wang}). However the proofs make 
no use of results from cyclotomic mappings.
Before giving a description of our work we need to set up  our notations.

\vspace{.5pc}
\noindent \textbf{Notation and Conventions:}

\vspace{1pc}
		\begin{tabular}{rll}
			\hline
		$q$    & :&     Power of an odd prime \strut  \\
		$\ff_q$ &:&  Finite field of order $q$ \\
		$g $ &:& a generator of $\ff_q^{\ast}$ \\
		$m,d$ &: &  pair of integers  factorizing $q-1$ as $q-1=md$\\
			$\varphi$&:& Euler totient function\\
		$H_0$ &:& the subgroup of $\ff_q^{\ast}$ of order $m$ (and hence index $d$)\\
		$H_i$ &:& coset of $H_0$ by $g ^{i}$, $i=0,1,\ldots,d-1$  \\ \hline
	     \end{tabular}

\vspace{1pc}
By the order of a (non-zero)  element, we always mean the order as an element of the \textit{multiplicative	group} of the field.      

Note that our definition of $m,d$ allows us to express  any element of $\ff_q^{\ast}$ 
 uniquely  as $g ^{i+dj}$ for $i=0,1,\ldots,d-1$ and $j=0,1,\ldots,m-1$.

 \vspace{6pt}
   With the notation set up  we are now  ready to give a description of our results.
For a finite field  $\ff_q,\ q$ odd, we construct  families of 
\begin{description}
\item[BIN-1]
binomial permutations  with zero as unique fixed-point, and 
remaining elements permuted as disjoint cycles of the same length. 

		For example for the prime field of 89 elements,
		we are able to construct 90 binomial permutations which
		permutate the nonzero elements into eight disjoint 
		cycles of length 11; another 90
		binomials representing  permutations of cycle type $1+22^4$.
		(see Remark 2 appearing after the proof of Theorem~\ref{thm1} and also
		the table at the end of the paper.)
\item[BIN-2] binomial
	permutations; for a fixed divisor $m$ of $q-1$, 
	and two  divisors $m_0, m_1$ of $m$ with $n_i= m/m_i$
		these binomials give permutations of cycle type
		$$
		1 + \underbrace{m_0+m_0+\cdots +m_0}_{dn_0/2 \mbox{ times}}
		 + \underbrace{m_1+m_1+\cdots +m_1}_{dn_1/2 \mbox{ times}}
		$$
	\item[TRI-1] for $q\equiv 1\pmod 3 $ trinomial permutations with zero as unique
		fixed-point and
remaining elements permuted as disjoint cycles of the same length. 
\item[TRI-2]  trinomial
	permutations; for a fixed divisor $m$ of $q-1$, 
	and three divisors $m_0, m_1, m_2$ of $m$ with $n_i= m/m_i$
		of cycle type
		$$
		1 + 
		  \underbrace{m_0+m_0+\cdots +m_0}_{dn_0/3 \mbox{ times}}
		+\underbrace{m_1+m_1+\cdots +m_1}_{dn_1/3 \mbox{ times}}
		 + \underbrace{m_2+m_2+\cdots +m_2}_{dn_2/3 \mbox{ times}}
		$$
	\item[POLY] and generalization of the above to polynomials of $n$ terms
		where $n$ is a divisor  of $q-1$.
\end{description}

The following elementary fact which is easy to prove will be repeatedly used:

\noindent \begin{lemma}\label{cyctype}
	Let $R$ be a commutative ring with unity. Assume the group of units of $R$
	has a finite subgroup $U$, say, of order $M$. Let $u\in U$, be an element of
	order $m$, and let $d=M/m$.
	Then, (i) the function $f\colon U\to U$ defined as $f(x) =ux$ is a permutation
	of $U$, (ii) $f(x)$ takes each coset of the cyclic subgroup generated by $u$ in $U$ to
	itself and (iii) the cycle type of the permutation given by $f(x) $ on $U$  is
	$$\underbrace{m+m+\cdots +m}_{d\ \rm terms}$$
\end{lemma}

Note:\quad In future we will use the simpler notation $m^d$ for the cycle type of the above permutation, and similar notation for general permutations.

\section{Main Results}

Our findings are spread over the four theorems below:
\begin{theorem}(\textbf{BIN-1})\label{thm1}
For  two distinct elements $u,v\in\ff_q^\ast$, 
	both of the same order $m$, with $d=(q-1)/m$ even,
	we define $f(x)\in \ff_q[x]$ as below:
	\[ f(x)=ax^{(q+1)/2}+bx\qquad \mbox{where}\quad 
	a=(v-u)/2,\ b=(u+v)/2.
	\]
	 Then the following are true for the  polynomial $f(x)$:
\begin{itemize}
	\item[(i)] it is a permutation binomial of $\ff_q$.
	\item[(ii)] the  permutation given by it has cycle type $1+m^d$
\item[(iii)] the inverse of this  permutation  is also a binomial given by 
	\[f^{-1}(x)=a'x^{\frac{q+1}{2}}+b'x  \]  where
	$a'=(v^{-1}-u^{-1})/2$ and $b'=(u^{-1}+v^{-1})/2$.
\end{itemize}
\end{theorem}
\begin{proof}
	We prove this directly by evaluating the polynomials, carrying out simple computations.
	\textit{We note for future use that 
\begin{equation}\label{eq1}
a+b=v ,\quad b-a=u
\end{equation}
 which follow immediately from their definitions. }
	For any $ j=0,1,2,\ldots,m-1$ we evalaute $f$ at $x=g^{i+dj}$ (note  that, as $d$
	is even, $i+dj$ is odd or even according as $i$ is odd or not.) 

	\vspace{6pt}
\noindent \textbf{Case  $i$ odd:}
	\begin{flalign*}
                f(g ^{i+dj})&=a(g ^{i+dj})^{\frac{q+1}{2}}+bg ^{i+dj}& \\
                &=g ^{i+dj}\bigl(a(g ^{i+dj})^{\frac{q-1}{2}}+b \bigr) &\\
		&=g ^{i+dj}(-a+b)\mbox{ (as odd powers of $g $ are non-squares)}&\\
                &=ug ^{i+dj} \mbox{ (by Eqn.\ \ref{eq1})} &
        \end{flalign*}

	\vspace{6pt}
	\noindent \textbf{ Case  $i$ even:}
	
	For any $ j=0,1,2,\ldots,m-1$ we have 
	\begin{flalign*}
		f(g ^{i+dj})&=a(g ^{i+dj})^{\frac{q+1}{2}}+bg ^{i+dj}& \\ 
        &=g ^{i+dj}\bigl(a(g ^{i+dj})^{\frac{q-1}{2}}+b \bigr) &\\
		&=g ^{i+dj}(a+b) \mbox{ (even powers of $g $ are squares)}&\\
		&=vg ^{i+dj}  \mbox{ (by Eqn.\ \ref{eq1})}
	\end{flalign*}
	This  shows that on the coset $H_i$ for $i$ odd, the binomial   $f(x)$ acts as multiplication by $u$ and  other cosets it acts as multiplication by $v$. That is
	\begin{equation} \label{eqn1}
		f(g ^{i+dj}) =
\left\{ \begin{array}{rcl}
   	ug ^{i+dj}  & \mbox{if  $i$ is  odd} \\
	vg ^{i+dj}  & \mbox{if  $i$ is even} \\

\end{array}\right.
	\end{equation}
As $u,v$ are elements of order $m$,  they  belong to $H_0$, the unique subgroup of order $m$.
	So $f(x)$ acts as full length cycle when restricted to each coset $H_i$. This shows
	$f(x)$ is a permutation polynomials, and from Lemma~\ref{cyctype}  its cycle decomposition
	is also easily deduced, proving
	parts (i) and (ii).

For part (iii),  replacing $u$ and $v$ in Eqn.\ (\ref{eqn1}) by  $u^{-1}$ and $v^{-1}$ we note that
we get the permutation inverse to $f(x)$. This is  also a binomial of the  same kind because in any group  an element and its inverse have the same order. 

\end{proof}

\noindent \textbf {Remark 1:}\quad  In \cite{nieder} authors construct binomials
where the coefficients are in similar format 
as $(v-u)/2$ and $(v+u)/2$. There the condition is simply that these $u,v$ 
be both squares or both non-squares. Here we have additional constraints that 
  both be of same order. But we reap the benefit of being able to
say that these constructions provide permutations of specific cycle type.

\noindent\textbf{Remark 2:}\quad This method shows, for any divisor $m$ of $q-1$, that there exist a family
of  (at least) $\varphi(m) (\varphi(m)-1)$
binomials of  degree  $(q+1)/2$ of cycle type $1+m^d$  when  $d=(q-1)/m$ is even. 

Now we want to state a variation where the two elements $u, v$ need not be of the same order.
\begin{theorem}\textbf{(BIN-2)}
Under the above notation let  $u,v$ be two different elements from the subgroup  $H_0$. Let their 
 orders be $m_0$ and $m_1$ respectively. Write $n_i=m/m_i,\ i=0,1$.
The binomial  \[ f(x)=ax^{\frac{q+1}{2}}+bx \] 
where $a=(v-u)/2$ and $b=(u+v)/2$
\begin{itemize}
\item[(i)] is a permution binomial on $\ff_q$.
\item[(ii)]the  permutation has cycle type $1+m_0^{dn_0/2}+m_1^{dn_1/2}$
\item[(iii)] the inverse of this  permutation  is also a binomial given by \[f^{-1}(x)=a'x^{\frac{q+1}{2}}+b'x  \]  where $a'=(v^{-1}-u^{-1})/2$ and $b'=(u^{-1}+v^{-1})/2$.

\end{itemize}

\end{theorem}
\begin{proof}
 The proof is similar to previous theorem with some modifications:  in  cosets $H_i$ for 
	$i$ odd, $f(x)$ is not a full cycle but a disjoint union of cycles having
	length $m_1$. Analagously  in  cosets $H_i$ for $i$ even, similar statement is 
	valid with $m_0$ replacing $m_1$. 
\end{proof}

Next we move on to construction of permutation binomials for fields $\ff_q$ in cases where 
$q-1$ is divisible by $3$. In previous cases we chose $d$ to be even; now we will
choose $d$ to be multiple of 3,
\begin{theorem}\textbf{(TRI-1)}
Let $q\equiv 1\pmod 3$ and $\zeta$ be an element of order 3
	in $\ff_q^*$. Consider a positive integer $d$ such that $d\equiv 0 \pmod 3$, and $q-1=md$ 
	with  $\varphi(m)\geq 3 $.
For any three elements $u,v,w\in\ff_q^\ast$, not all same,   chosen to be of same order $m$, the trinomial  \[ f(x)=ax^{\frac{2q+1}{3}}+bx^{\frac{q+2}{3}}+cx \] 
	where $a=(u+v\zeta+w\zeta^2)/3$, $b=(u+v\zeta^2+w\zeta)/3$ and
	$c=(u+v+w)/3$
\begin{itemize}
\item[(i)] is a permutation polynomial on $\ff_q$.
\item[(ii)]this  permutation has cycle type $1+ m^d$.
\item[(iii)] its inverse permutation  is also a trinomial given by \[f^{-1}(x)= a'x^{\frac{2q+1}{3}}+b'x^{\frac{q+2}{3}}+c'x \]
where $a'=(u^{-1}+v^{-1}\zeta+w^{-1}\zeta^2)/3$ , $b'=(u^{-1}+v^{-1}\zeta^2+w^{-1}\zeta)/3$ and $c'=(u^{-1}+v^{-1}+w^{-1})/3$.

\end{itemize}
\end{theorem}
\begin{proof}
As $d\equiv0\pmod 3$ we have  $\zeta^{d}=1$.
To  prove $f(x)$ is  a permutation polynomials we need to discuss the following cases with respect to the exponent of $g $.

\noindent 
	Let us evaluate the polynomial $f(x)$ for $x$ in $\ff_q^*$ expressing it in the form
	 $x = g^{i+dj}$. As $d$ is a multiple of 3, $i+dj$ and $i$ will be same modulo 3.

	\vspace{6pt}
\noindent\textbf{Case  $i\equiv 0\pmod 3$ }

In this case we write $i=3k$. 
	Now $f(x)$ at $g ^{3k+dj}$ is
\begin{flalign*}
f(g ^{3k+dj})&=g ^{3k+dj}\bigl(a(g ^{3k+dj})^{ \frac{2(q-1)}{3} }+b(g ^{3k+dj})^{ \frac{q-1}{3} } +c \bigr) \\
                &=g ^{3k+dj}\bigl(a(g ^{dj})^{\frac{2(q-1)}{3}}+b(g ^{dj})^{\frac{q-1}{3}} +c \bigr) \\
                &=g ^{3k+dj}\bigl(a\zeta^{2dj}+b\zeta^{dj} +c \bigr) \\
                &=g ^{3k+dj}\bigl(a+b +c \bigr) \\
                &=ug ^{3k+dj}
\end{flalign*}
where the last equality is obtained by substituting the values of $a$, $b$ and $c$ .

	\vspace{6pt}
 \noindent \textbf{Case  $i\equiv 1\pmod 3$ }	

In this case we write $i=3k+1$.
Now $f(x)$ at $x=g ^{3k+1+dj}$ is
\begin{flalign*}
f(g ^{3k+1+dj})&=g ^{3k+1+dj}\bigl(a(g ^{3k+1+dj})^{ \frac{2(q-1)}{3} }+b(g ^{3k+1+dj})^{ \frac{q-1}{3} } +c \bigr) \\
                &=g ^{3k+1+dj}\bigl(a(g ^{1+dj})^{\frac{2(q-1)}{3}}+b(g ^{1+dj})^{\frac{q-1}{3}} +c \bigr) \\
                &=g ^{3k+1+dj}\bigl(a\zeta^{2dj}\zeta^2+b\zeta^{dj}\zeta +c \bigr) \\
                &=g ^{3k+1+dj}\bigl(a\zeta^2+b\zeta +c \bigr) \\
                &=vg ^{3k+1+dj}
\end{flalign*}
here again the last equality is obtained by substituting the values of $a$,$b$ and $c$.

\vspace{6pt}
\noindent \textbf{Case  $i\equiv 2\pmod 3$ }

In this case we write $i=3k+2$. 
	Now $f(x)$ at $x=g ^{3k+2+dj}$ is
\begin{flalign*}
f(g ^{3k+2+dj})&=g ^{3k+1+dj}\bigl(a(g ^{3k+2+dj})^{ \frac{2(q-1)}{3} }+b(g ^{3k+2+dj})^{ \frac{q-1}{3} } +c \bigr) \\
                &=g ^{3k+2+dj}\bigl(a(g ^{2+dj})^{\frac{2(q-1)}{3}}+b(g ^{2+dj})^{\frac{q-1}{3}} +c \bigr) \\
                &=g ^{3k+2+dj}\bigl(a\zeta^{2dj}\zeta+b\zeta^{dj}\zeta^2 +c \bigr) \\
                &=g ^{3k+2+dj}\bigl(a\zeta+b\zeta^2 +c \bigr) \\
                &=wg ^{3k+2+dj}
\end{flalign*}
	So the polynomial $f(x)$ represents the  function admitting the description below:

	\begin{equation} \label{eqn2}
		f(g ^{i+dj}) =
\left\{ \begin{array}{rcl}

	ug ^{i+dj}  & \mbox{ if   $i \equiv 0 \pmod 3$}  \\
	vg ^{i+dj}  & \mbox{ if   $i \equiv 1 \pmod 3$} \\
	wg ^{i+dj}  & \mbox{ if   $i \equiv 2 \pmod 3$} \\
\end{array}\right.
	\end{equation}
Now it is clear that  $f(x)$ is a permutation polynomials and its cycle type can again deduced
from Lemma~\ref{cyctype} to be 
$1+m^d$.
	Replacing $u, v$ and $w$ in Eqn.\  (\ref{eqn2}) by  $u^{-1}, v^{-1}$ and $w^{-1}$ we get the
	trinomial permutation  which is the  inverse of $f(x)$.

\end{proof}

\begin{theorem}\textbf{(TRI-2)}
Let $q\equiv 1 \pmod 3$ and $\zeta$ be an element of order 3. Consider the positive integers $m,d $ such that $d\equiv 0 \pmod 3$, $q-1=md$.
For any three elements $u,v,w\in\ff_q^\ast$ whose orders are $m_0$,$m_1$ and $m_2$ respectively. Write $n_i=\frac{m}{m_i},\ i=0,1,2$.

The  trinomial  \[ f(x)=ax^{\frac{2q+1}{3}}+bx^{\frac{q+2}{3}}+cx \] 
where $a=(u+v\zeta+w\zeta^2)/3$, $b=(u+v\zeta^2+w\zeta)/3$ and
	$c=(u+v+w)/3$
\begin{itemize}
\item[(i)] is a permution trinomial on $\ff_q$
\item[(ii)]the  permutation has cycle type \[ 1+ m_0^{dn_0/3}
	+ m^{dn_1/3}
		+ m_2^{dn_2/3}
		\] 
\item[(iii)] its inverse permutation also trinomial given by \[f^{-1}(x)= a'x^{\frac{2q+1}{3}}+b'x^{\frac{q+2}{3}}+c'x \]
where $a'=(u^{-1}+v^{-1}\zeta+w^{-1}\zeta^2)/3$ , $b'=(u^{-1}+v^{-1}\zeta^2+w^{-1}\zeta)/3$ and $c'=(u^{-1}+v^{-1}+w^{-1})/3$.

\end{itemize}
\end{theorem}
\begin{proof}

Proof is similar to previous theorem.
\end{proof}
\section{Conclusion}
The results proved above have generalization to constructing polynomials
having  $r$ terms for any $r$ dividing $q-1$ with predictable cycle type.
The coefficients of such polynomials
are given by uncomplicated formula coming from  the solutions of a 
linear system corresponding to  an easy special case of  Vandermonde  matrix.

\begin{theorem}\textbf{(POLY)} \label{thrm1}
Let $r>1$ be a positive integer. Consider the finite field $\ff_q$ with
$q\equiv 1 \pmod r$ and let $\zeta$ be an element of order  $r$. We choose $m,d$
	satisfying   
	$q-1=md$,  $d$ be a multiple of $r$ and $\varphi(m)\geq r$. 
	For any $r$ elements $u_0,u_1,\ldots,u_{r-1}$ from $\ff_q^\ast$ all of the same  order $m$,  we define a  polynomial
\[G(x)=a_{r-1}x^{r-1}+a_{r-2}x^{r-2}+\cdots+a_1x+a_0\in\ff_q[x] \] whose coefficients are
	determined from the  following system of  equations: 
\begin{equation*}\label{eqn3}
\left(\begin{array}{rrrrrrr}
1                   &1                & \cdots   &1          &1          &1 \\
 \zeta^{r-1}       &\zeta^{r-2}     &\cdots    &\zeta^2   &\zeta     & 1\\
\zeta^{2(r-1)}     &\zeta^{2(r-2)}  &\cdots    &\zeta^4   &\zeta^2   &1\\
\vdots              &\vdots           &\cdots    &\vdots     & \vdots    & \vdots \\
\zeta^{(r-1)(r-2)} &\zeta^{(r-2)(r-2)} &\cdots &\zeta^{2(r-2)}  &\zeta^{r-2} &1 \\
\zeta^{(r-1)(r-1)} &\zeta^{(r-1)(r-2)} &\cdots &\zeta^{2(r-1)}  &\zeta^{r-1} &1
\end{array}\right) 
\left(\begin{array}{rrrrrrr}
a_{r-1} \\ a_{r-2} \\ a_{r-3}\\ \vdots \\ a_1  \\ a_0
\end{array}\right) 
=
\left(\begin{array}{rrrrrr}                 
u_0 \\ u_1 \\u_2 \\ \vdots \\ u_{r-2} \\ u_{r-1} \end{array}\right) 
\end{equation*}

	\begin{equation} \label{vander}
\left(\begin{array}{rrrrrrr}
a_{r-1} \\ a_{r-2} \\ a_{r-3}\\ \vdots \\ a_1  \\ a_0
\end{array}\right) 
=
\left(\begin{array}{rrrrrrr}
1        &\zeta     &\zeta^2   & \cdots   &\zeta^{r-2}          &\zeta^{r-1} \\
 1       &\zeta^2   &\zeta^4   &\cdots    &\zeta^{2(r-2)}   &\zeta^{2(r-1)}   \\
 1       &\zeta^3   &\zeta^6   &\cdots    &\zeta^{3(r-2)}   &\zeta^{3(r-1)}   \\

\vdots   &\vdots   &\vdots         &\cdots    &\vdots     & \vdots  \\
 1       &\zeta^{r-1}   &\zeta^{2(r-1)}   &\cdots    &\zeta^{(r-1)(r-2)}   &\zeta^{(r-1)(r-1)}   \\
1 &  1 & 1 &\cdots &1  &1 
\end{array}\right) 
\left(\begin{array}{rrrrrrr}                 
u_0 \\ u_1 \\u_2 \\ \vdots \\ u_{r-2} \\ u_{r-1} \end{array}\right)   
\end{equation}

Then $f(x)=xG(x^{\ell})$ where $\ell=\frac{q-1}{r}$,
\begin{itemize}
\item[(i)]  is a permutation polynomial on $\ff_{q}$.
\item[(ii)] the  permutation has cycle type $1+m^d$.
\item[(iii)]  its inverse is also a polynomial of $n$ terms involving the same powers of $x$ as $f(x)$. Its coefficients are determined 
	from  the Vandermonde system (\ref{vander}) by replacing $u_i, 0\leq i\leq r-1 $  with its inverse
\end{itemize}
\end{theorem}
\begin{proof}
 We can easily show the polynomial $f(x)$ is permutation polynomials from its function respresentation  below:
\[	f(g ^{i+dj}) = u_{i}g ^{i+dj}\] here $0\leq i \leq r-1$.
\end{proof}

Again we can state a variation of the above theorem using the same notations:
\begin{theorem}  
 For any $n$ elements $u_i$ from $\ff_q^\ast$, not all the same,  with $u_i$ of order
	$m_i$ dividing $m$, write $n_i=\frac{m}{m_i},\  0\leq i\leq r-1$.
 The polynomial $f(x)$  given in Theorem (\ref{thrm1}).

\begin{itemize}
\item[(i)]  is a permutation polynomial on $\ff_{q}$.
\item[(ii)] the  permutation has cycle type
	\[1+ 	m_0^{dn_0/n} +
		m_1^{dn_1/n} +\cdots+
		m_n^{dn_i/n} \]
\item[(iii)]  its inverse permutation is also a polynomial 
	of $n$ terms involving the same powers of $x$ as $f(x)$. Its coefficients are determined from the system (\ref{vander}) modified 
		by replacing  $u_i$'s with their inverse.
\end{itemize}
\end{theorem}

Using the same convention for notations explained in the introduction we have the following
result:
\begin{theorem}
	For two distinct elements
	$u,v\in\ff_q^\ast$, both of the same order $m$ define $a= \frac{u-v}{d}$. Then the
	polynomial \[f(x)= a\sum_{i=0}^{d-1} x^{im+1}+vx\] is  
	a permutation polynomial giving a permutation of cycle type $1+m^d$. 
\end{theorem}

\begin{proof}
Since given polynomial takes the value zero at zero, enough to show $f(x)$ permutes all non-zero 
elements which is shown by considering two  cases:
Let $b\in\ff_q^\ast$.

\vspace{6pt}
\noindent \textbf{Case:  $ b\in H_0$ }
\begin{flalign*}
f(b) &= a(b+b^{m+1}+b^{2m+1}+\cdots+b^{(d-1)m+1})+vb &\\
     &= ab(1+b^m+\cdots+b^{(d-1)m})+vb &\\
	&=abd+vb \mbox{ [since $b^m=1$]} &\\
	&=ub \mbox{ [substituting  value of $a$.] } &\\ 
\end{flalign*}

\vspace{3pt}
\noindent \textbf{Case:  $ b\notin H_0$ }
\begin{flalign*}
f(b) &= a(b+b^{m+1}+b^{2m+1}+\cdots+b^{(d-1)m+1})+vb &\\
     &= ab(1+b^m+\cdots+b^{(d-1)m})+vb &\\
     &= vb \mbox{ (since  $1+b^m+\cdots+b^{(d-1)m}=0$)}  &\\
\end{flalign*}
 this completes the proof for the fact that  $f(x)$ is a permutation 
	polynomial;  the cycle type is deduced by again appealing to
	Lemma~\ref{cyctype}. The statement about the compositonal inverse  of this polynomial
	is also proved by arguments used earlier.

\end{proof}

\noindent \textbf{Example:}
We tabulate below, for the finite field of $121=11^2$ elements,
the count of binomials and trinomials we were able to construct
having specific cycle type.
Each row of the table corresponds to a choice of factorization of 120= $m\times d$
leading to a specific conjugacy class of permutations.

\begin{center}
\begin{tabular}{|l|r|r|c|r|}
\hline
	\textbf{Polynomial type}
	& $m$ &  $d$ \vrule height13pt depth5pt width0pt
	& \textbf{Cycle type}  & \textbf{Count} \\ \hline 
\multirow{10}{*}{Binomial}  & 3 & 40 & $1+3^{40}$ & 2 \vrule height13pt depth5pt width0pt\\ \cline{2-5}
                            & 4 & 30 & $1+4^{30}$ & 2 \vrule height13pt depth5pt width0pt\\ \cline{2-5}
                            & 5 & 24 & $1+5^{24}$ & 12 \vrule height13pt depth5pt width0pt\\ \cline{2-5}
              				& 6 & 20 & $1+6^{20}$ & 2 \vrule height13pt depth5pt width0pt\\  \cline{2-5}
							& 10 & 12 & $1+10^{12}$ & 12 \vrule height13pt depth5pt width0pt\\ \cline{2-5} 
						& 12 & 10 & $1+12^{10}$ & 12 \vrule height13pt depth5pt width0pt\\  \cline{2-5}
 							& 15 & 8 & $1+15^{8}$ & 56 \vrule height13pt depth5pt width0pt\\  \cline{2-5}
 							& 20 & 6 & $1+20^{6}$ & 56\vrule height13pt depth5pt width0pt \\  \cline{2-5}
 							& 30 & 4 & $1+30^{4}$ & 56\vrule height13pt depth5pt width0pt \\  \cline{2-5}
							& 60 & 2 & $1+60^2$   & 240\vrule height13pt depth5pt width0pt \\ 
\hline

\multirow{5}{*}{Trinomial} 
& 5 & 24 & $1+5^{24}$ & 60 \vrule height13pt depth5pt width0pt\\ \cline{2-5}
& 10 & 12 & $1+10^{12}$ & 60 \vrule height13pt depth5pt width0pt\\ \cline{2-5} 
& 8 & 15 & $1+8^{15}$ & 60 \vrule height13pt depth5pt width0pt\\  \cline{2-5}
& 20 & 6 & $1+20^{6}$ & 504\vrule height13pt depth5pt width0pt \\ \cline{2-5}
& 40 & 3 & $1+40^{3}$ & 4080\vrule height13pt depth5pt width0pt \\
 
\hline   
\end{tabular}
\end{center}


\end{document}